\documentclass[11pt]{article}
	
	\newcommand{\blind}{0}
	
    \makeatletter
    \renewcommand\section{\@startsection {section}{1}{\z@}%
                                       {-3.5ex \@plus -1ex \@minus -.2ex}%
                                       {2.3ex \@plus.2ex}%
                                       {\normalfont\fontfamily{phv}\fontsize{16}{19}\bfseries}}
    \renewcommand\subsection{\@startsection{subsection}{2}{\z@}%
                                         {-3.25ex\@plus -1ex \@minus -.2ex}%
                                         {1.5ex \@plus .2ex}%
                                         {\normalfont\fontfamily{phv}\fontsize{14}{17}\bfseries}}
    \renewcommand\subsubsection{\@startsection{subsubsection}{3}{\z@}%
                                        {-3.25ex\@plus -1ex \@minus -.2ex}%
                                         {1.5ex \@plus .2ex}%
                                         {\normalfont\normalsize\fontfamily{phv}\fontsize{14}{17}\selectfont}}
    
    \makeatother

	\usepackage{amsmath}
	\usepackage{graphicx}
	\usepackage{enumerate}
	\usepackage{xcolor}
    \usepackage{subcaption}
    \usepackage{hyperref}
    \usepackage{fancyhdr}
	\usepackage{natbib} 
	\usepackage{url} 

	\DeclareMathOperator{\prob}{Pr}
	
\begin{document}		
		\def\spacingset#1{\renewcommand{\baselinestretch}%
			{#1}\small\normalsize} \spacingset{1}
		
		\if0\blind
		{
			\title{\bf Co-optimization of Short- and Long-term Decisions for the Transmission Grid's Resilience to Flooding}
			\author{Ashutosh Shukla$^{\dagger}$ , Erhan Kutanoglu, John Hasenbein \\ Graduate Program in Operations Research and Industrial Engineering, \\The University of Texas at Austin, Austin, Texas, USA}
			\date{}
			\maketitle
		} \fi
		
		\if1\blind
		{
            \title{\bf Co-optimization of Short- and Long-term Decisions for the Transmission Grid's Resilience to Flooding}
			\author{}
			\maketitle    
		} \fi
		\bigskip		
	\begin{abstract}
We present and analyze a three-stage stochastic optimization model that integrates output from a geoscience-based flood model with a power flow model for transmission grid resilience planning against flooding. The proposed model coordinates the decisions made across multiple stages of resilience planning and recommends an optimal allocation of the overall resilience investment budget across short- and long-term measures. While doing so, the model balances the cost of investment in both short- and long-term measures against the cost of load shed that results from unmitigated flooding forcing grid components go out-of-service. We also present a case study for the Texas Gulf Coast region to demonstrate how the proposed model can provide insights into various grid resilience questions. Specifically, we demonstrate that for a comprehensive yet reasonable range of economic values assigned to load loss, we should make significant investments in the permanent hardening of substations such that we achieve near-zero load shed. We also show that not accounting for short-term measures while making decisions about long-term measures can lead to significant overspending. Furthermore, we demonstrate that a technological development enabling to protect substations on short notice before imminent hurricanes could vastly influence and reduce the total investment budget that would otherwise be allocated for more expensive substation hardening. Lastly, we also show that for a wide range of values associated with the cost of mitigative long-term measures, the proportion allocated to such measures dominates the overall resilience spending.
	\end{abstract}
			
	\noindent%

	\spacingset{1.5} 

\vfill
\hrule
\vspace{1em}
\noindent $\dagger$ Corresponding author, now employed at BNSF Railway, Fort Worth, Texas, USA\\
\newpage
\section{Introduction}

\subsection{Motivation}

Tropical cyclones, severe storms, and inland flooding have collectively inflicted over \$2 trillion (CPI-adjusted) in damages between 1980 and 2023 in the U.S \citep{billion_nerc}. More than half of these damages come from the five states (Texas, Louisiana, Mississippi, Alabama, and Florida) exposed to the Gulf of Mexico. A significant portion of these losses can be attributed to the damage of critical infrastructures such as the power grid. For example, Hurricane Harvey, which affected the Texas coast between Corpus Christi and the Houston-Galveston region in 2017, damaged more than 90 substations, 800 transmission structures, 6000 distribution poles, and 800 miles of transmission and distribution lines. The storm caused a peak generation loss of 11 GW and affected more than two million customers. It took two weeks and 12,000 crew members to restore power \citep{NERC2018}. The failure of the power grid impacted the functioning of other life-line infrastructures like healthcare services, transportation networks, critical supply chains, and water networks, leading to significant social and economic losses.

Given the importance of protecting the power grid from extreme weather events, several studies have focused on developing methodological frameworks \citep{framework_grid_resilience, Panteli_framework, framework_5} and metrics \citep{Panteli_metrics} for power grid resilience to extreme weather events. Additionally, various stochastic programming-based models have been developed that combine scenarios representing the uncertainty of damage to power grid components with various representations of power flow models for resilience decision-making. These models typically weigh the upfront investment cost to mitigate the impact of hurricanes against the social and economic costs incurred due to the curtailed performance of the grid with insufficient protection. For example, the two-stage stochastic optimization model proposed in \cite{7540988} considers the trade-off between the upfront costs associated with the upgrade of the transmission system by hardening existing components and adding redundant lines, switches, generators, and transformers in the first-stage against the benefits of these upgrades provide in the second stage against flooding. Similarly, the model proposed in \cite{9376648} suggests an investment plan to minimize the total investment cost of building new transmission lines and deployment of distributed energy resources, expected values of generators’ operation cost in normal and emergency situations, and load shedding cost in emergency conditions. Another model proposed in \cite{9470975} recommends an optimal plan for hardening of transmission lines, generators, and substations such that the load loss during the recovery phase is minimized. The two-stage model proposed in \cite{9842364} also focuses on line hardening decisions but accounts for the impact of decision-dependent uncertainty associated with typhoon disasters. In addition to grid hardening, the resilience of the transmission system can be further enhanced through grid islanding techniques as demonstrated in \cite{Panteli_boosting}. Similar to the stochastic programming models, more recently, a two-stage robust optimization model based on the $N-k$ security criterion has been proposed in \cite{ZHANG2024110068} for power grid resilience decision-making.

In generating the scenarios to represent the impact of hurricanes, the aforementioned models do not consider the impact of flooding. However, it is imperative to recognize that hurricanes have the potential to inflict substantial damage on substations through flooding. Towards that end, the models proposed in \citep{Mohadese, austgen_impacts_2021, 10252826, SOUTO2022107545, shukla_scenario-based_2021} focus explicitly on capturing the impact of flooding on the power grid and recommend strategies for substation protection. The models in \cite{Mohadese} and \cite{10252826} use fragility curves based on log-normal distributions to represent the risk of substation failures, whereas the models proposed in \cite{SOUTO2022107545} and \cite{austgen_impacts_2021} use a physics-based simulation model to predict the flood risk for simulated storm scenarios. Specifically, \cite{Mohadese, 10252826, austgen_impacts_2021} focus on decision-making in the preparedness phase for an impending but still uncertain storm and recommend a substation protection strategy using ad hoc, temporary measures such as Tiger Dams\textsuperscript{\texttrademark} (the superscript is omitted in further discussion for a clearer presentation) shown in Figure \ref{fig:Tiger_dams}. On the other hand, the models proposed in \cite{SOUTO2022107545} and \cite{shukla_scenario-based_2021} focus on longer-term permanent substation hardening measures against flooding due to multiple future storms. However, while making such hardening decisions, these models do not account for preventive measures such as Tiger Dams that can be deployed before an imminent hurricane's landfall.

\subsection{Research Gap}

The lack of coordination between the short- and long-term resilience decisions is the research gap that this study aims to address. The benefits of coordinated decision-making between the different phases of resilience management has been demonstrated in a limited number of studies. For example, in \cite{HUANG2024110136}, the authors propose a coordination of preventive, emergency and restoration dispatch (CPERD) methodology to counteract cascading failures such that both the total load loss and dispatch costs are minimized. In an another study, a two-stage robust optimization model is proposed in \cite{7885130} that integrates measures like generation dispatch, transmission switching, and load-shedding, and shows that integrated resilience planning is preferable over doing preparedness and response planning independently. Similarly, the impact of coordinating preventive, emergency, and restorative power dispatch in transmission systems is shown using a hierarchical defender-attacker model in \cite{9591362}. Both are robust optimization models and consider only line failures.

To the best of our knowledge, no study exists that demonstrates the advantage of coordinated resilience decision-making to prevent substation failures due to hurricane-induced flooding. Motivated by these observations, we present a three-stage stochastic programming model for informing investment decision-making for the transmission grid, with a specific focus on substation failures. The flood uncertainty in the model is represented using scenarios generated from the output of a physics-based storm-surge model. The majority of the literature and most of the cited work discussed so far rely on statistical models to represent weather or flood uncertainty, and eventually scenarios. The exception to this are studies in \cite{SOUTO2022107545, shukla_scenario-based_2021, austgen_impacts_2021}, which integrate forecasts from weather or flood models as scenarios, with models of power flow. Even then, in contrast to these related studies, our model jointly optimizes long- and short-term investments for transmission grid resilience planning against substation failures caused by flooding. 

\subsection{Contributions}
In this study, we make the following contributions:
\begin{enumerate}
    \item We propose a three-stage stochastic optimization model that co-optimizes permanent substation hardening decisions made in the mitigation phase with Tiger Dam deployment decisions in the preparedness phase. Solving this model allows, among other things, to determine an optimal budget allocation between long- and short-term measures to enhance the transmission grid resilience. While doing so, the model's third (recourse) phase accounts for economic costs associated with load loss due to unmitigated flooding incurred during the planning horizon.
    \item We capture the impact of correlated flooding on the grid using a flood impact-sensitive DC power flow model that is embedded within the three-stage model.
    \item We further develop a case study for the state of Texas using a 663-bus synthetic grid that represents the transmission grid in the Texas Gulf Coast. Similarly, to represent the storm-surge-induced flood-risk on the coastal grid, we use the flood scenarios developed by the National Oceanic and Atmospheric Administration (NOAA). 
    \item We use the proposed model with the case study to answer a wide variety of resilience questions with an extensive sensitivity analysis on the value of load loss, power restoration time, ad hoc flood barrier's effectiveness, and substation hardening cost. We characterize the effects of these factors on the optimal allocation of the resilience budget across substation hardening, Tiger Dam acquisition, and deployment as a function of the value of load loss.
\end{enumerate}

\section{Modeling}
\subsection{Overview}
The proposed three-stage stochastic optimization model coordinates the decisions made in the different stages of resilience management to minimize the disaster cost of flooding of the transmission grid components. These decisions concern the measures that can be taken at different time scales in the resilience management cycle. We consider two such measures. The first is the permanent hardening of the substations, which is done to mitigate the effects of not one specific storm but of many potential storms that might affect the power grid during the multi-year planning horizon. The second is the deployment of ad hoc flood barriers as a preventative measure before an imminent hurricane's landfall. 

\begin{figure}
    \centering
    \includegraphics[scale=0.4]{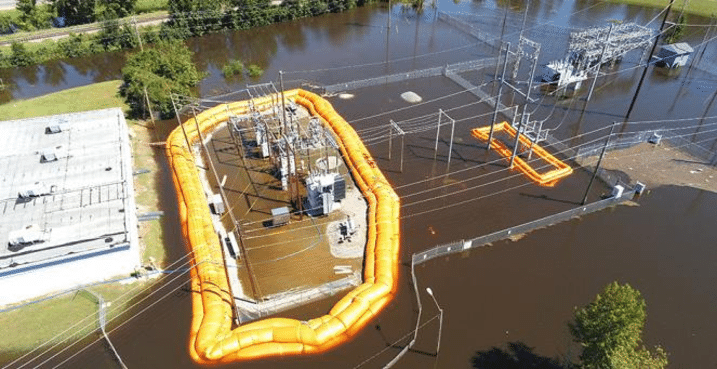}
    \caption{Tiger Dam protecting a substation in North Carolina after Hurricane Florence}
    \label{fig:Tiger_dams}
\end{figure}

Examples of permanent hardening measures include elevation of certain substation components and construction of permanent walls.
In either case, the resilience level of a substation can be expressed in terms of the height of flooding that it can withstand while continuing to remain operational. In the sequel, a resilience level of $l$ denotes the ability of a substation to withstand $l$ feet of flooding. Ad hoc barriers such as Tiger Dams (shown in Figure \ref{fig:Tiger_dams}) perform a similar function as that of permanent hardening but offer additional flexibility because they can be deployed quickly at a safe time before a landfall to temporarily protect substations subject to flooding threat due to the impending storm. This helps us account for the most-up-to-date data on the hurricane characteristics, potentially with more certainty as compared to the time of permanent hardening which is done months or years before any specific hurricane. Tiger Dams, however, have two limitations as compared to permanent hardening. The first is that they need to be disassembled and redeployed before every storm because they undergo wear and tear due to environmental factors if left unattended after a hurricane. As a consequence, the repeated disassembly and deployment of the dams for every hurricane incurs a recurring deployment cost. The second limitation is that Tiger Dams have a threshold lower than that of permanent hardening on the level of resilience they  provide. We now specify the types of decisions we consider in our model and the associated timeline. 

\textbf{First-stage decisions:} The first-stage mitigation phase decisions in the proposed model are the permanent hardening of the substations and the acquisition of the Tiger Dams. Specifically, the model decides (1) which substations to permanently harden and to what height to protect them against flooding, and (2) the number of Tiger Dam units to purchase to be used for deployment in the second stage before individual hurricanes. These first-stage decisions are made at the beginning of the multi-year planning horizon well before any specific hurricane. These decisions attempt to minimize the total expected disaster cost incurred during the planning horizon due to the flooding of the power grid components by a number of hurricanes with varying characteristics captured through a set of scenarios. 

\textbf{Second-stage decisions:} Representing the preparedness phase before an impending storm, the second stage decisions concern the deployment of Tiger Dams purchased in the first stage. During this stage, more information about a specific storm is assumed to be available and the uncertainty about the storm characteristics is reduced. As a consequence, the model considers a reduced set of scenarios as compared to the first-stage to decide where the Tiger Dams should be deployed and to what height. 

\textbf{Third-stage decisions:} In the third and final stage, representing the response phase, one of the scenarios is realized, causing a subset of substations with insufficient protection to go out-of-service due to flooding. As a consequence, the buses and the branches connected to the flooded substations are also out-of-service, affecting the capability of the grid to satisfy loads. Consequently, the last set of decisions in the recourse stage are about the power flow through the remaining operational grid with the goal of satisfying the loads as much as possible, or equivalently, minimizing the load loss.

In order to better understand our stochastic optimization formulation, we discuss a few assumptions regarding storms that occur over the planning horizon. We assume that an unknown (random) number of storms occur during the horizon. In addition, each storm occurrence can be thought of as being selected independently from the set of storm scenarios described in Section \ref{section_scenario_tree}. Hence, our assumption is that the storm distribution (i.e., the probability mass function over the set of scenarios) is time homogeneous. A more sophisticated model could attempt to incorporate climate change or other global effects that alter storm characteristics over time. However, we believe that, in general, the variability in the intensity and location of individual storms is much greater than the induced temporal variability, at least for the relatively short horizon we consider. Due to our assumptions, the objective function is unaffected by the timing or order of the storms. In addition, in computing average costs due to power outages, note that we can simply compute the cost averaged over all storm scenarios (via third stage), and then multiply by the average number of storms expected during the planning horizon, due to linearity of the expectation operator. Accordingly, we assume that on average, $Y$ hurricanes make landfall during the planning horizon, and before each landfall, we deploy Tiger dams considering the most up-to-date data on the individual hurricanes. These hurricanes can make landfall from different directions, have different forward speeds, and intensity categories. Furthermore, for each hurricane type (which is uniquely defined by the direction of landfall, intensity category, and forward  speed), we have an associated flood map which, when overlayed with the physical grid data, quantifies the flood heights that the substations face.

We distinguish between the uncertainty faced while making the first- and the second-stage decisions and represent them using random variables $\xi_1$ and $\xi_2$, respectively. Specifically, while making the first-stage decisions (about substation hardening and tiger dam acquisition) at the beginning of the planning horizon, we consider a large set of hurricanes with different characteristics that are likely to make landfall during the planning horizon.  This set represents the sample space of $\xi_1$. In the second-stage, before making the tiger dam deployment decisions, typically a few days ahead of when the hurricane is anticipated to make a landfall, we assume that the uncertainty about the hurricane characteristics is greatly reduced. Accordingly, the sample space of $\xi_2$ is similar to that of $\xi_1$ but considers much fewer possible hurricanes. Moreover, to determine the expected cost of tiger dam deployment and the economic cost of load shedding, we compute the average of these values for different realizations of $\xi_1$ and $\xi_2$ and then multiply it with $Y$ to represent the average cost for $Y$ hurricanes.

Given that decision makers typically have limited funds for mitigation measures like substation hardening and Tiger Dams which can be deployed at various levels providing different levels of protection, it is reasonable to trade off their costs and effectiveness against the total value of load loss. For this purpose, a unit value is assigned to the load loss to represent the incurred social cost due to power outages. An important question then is how much should be invested in permanent hardening and in Tiger Dams so that the total costs, including costs associated with the deployment of Tiger Dams and social costs due to load loss, are minimized. Our proposed model is designed to explore this trade-off under different cost structures. For a given set of parameters, it recommends an optimal allocation of funds for permanent hardening and purchase of Tiger Dams so as to minimize the aforementioned cumulative sum of costs across permanent hardening, preparedness, and recovery stages. 

\subsection{Assumptions}
To ensure that the proposed model can be used for studying large grid instances while considering a wide variety of flooding scenarios, we make a few assumptions to maintain computational tractability. The first assumption is that a DC approximation of the power flow as described in \cite{dc_opf_iie} is acceptable when embedding it within a larger power grid resilience decision-making model. This approximation has been embedded within decision-making models for resiliency-oriented transmission planning in several prior studies \cite{9376648, 10287582}. The second assumption is that the Tiger Dams units can be re-used during the planning horizon and do not undergo any degradation in performance. Furthermore, the cost incurred due to their deployment remains constant during the planning horizon. This is reasonable since they come in distinct units, and a dam is constructed by joining such units and stacking them on top of each other. Each of these individual units is typically filled with water. In the aftermath of the hurricane, the water is discharged, and the dams are disassembled into smaller units and stored until they are needed again. Third, we assume that each substation's hardening, Tiger Dam, and flood levels are discrete and finite. In the proposed model, we consider them to be non-negative integers. Fourth, we assume the grid experiences load loss due to substation flooding until it fully recovers after a storm.  Fifth, we assume that the value of load loss can be quantified in dollars per hour and remains constant during the planning horizon. The last two assumptions allow us to calculate the total load loss and its associated cost without explicitly modeling the detailed recovery schedule of the grid. 

\subsection{Notation}
In the following notation, all the cost parameters are in dollars. The power grid parameters associated with the branch flow capacity, power demand at a bus, and minimum and maximum power generation at buses are in the per-unit system.\\
\noindent\textbf{Random Variables}
\begin{itemize}
    \item[] $\xi_1, \xi_2$: Discrete random variables associated with stage one and stage two uncertainty, respectively 
\end{itemize}
\noindent\textbf{Sets}
\begin{itemize}
    \item[] $\mathcal{I}_f$: Set of flooded substations
    \item[] $\mathcal{H}_i, \mathcal{T}_i$: Set of integer hardening and Tiger Dam levels at substation $i$
    \item[] $\mathcal{J}$: Set of buses indexed by $j$
    \item[] $\mathcal{B}_i$: Set of buses at substation $i$
    \item[] $\mathcal{R}$: Set of branches indexed by $r$
    \item[] $\mathcal{N}_j^{in}$, $\mathcal{N}_j^{out}$: Set of branches incident on bus $j$ with power flowing in and out of bus $j$, respectively
    \item[] $\mathcal{P}$: Sample space of $\xi_1$
    \item[] $\mathcal{Q}_k$: Sample space of $\xi_2$ when $\xi_1 = k$
\end{itemize}
\noindent\textbf{Parameters}
\begin{itemize}
    \item[] $M$: An arbitrarily large constant
    \item[] $V_{il}$: Cost of hardening substation $i$ to resilience level $l$
    \item[] $A$: Acquisition cost of a Tiger Dam unit
    \item[] $C_{il}$: Cost of deploying a Tiger Dam at substation $i$ to prevent up to $l$ feet of flooding
    \item[] $Y$: Average number of storms during the planning horizon
    \item[] $W_i$: Maximum flood height at substation $i$ across all the scenarios
    \item[] $\Delta_{ikm}$: Flood height at  substation $i$ when $\xi_1 = k$, $\xi_2 = m$
    \item[] $B_r$: Susceptance of branch $r$
    \item[] $P_r$: Maximum power that can flow in branch $r$
    \item[] $r.\mathrm{h}$, $r.\mathrm{t}$: Head bus and tail bus of branch $r$
    \item[] $D_j$: Power demand at bus $j$
    \item[] $\underline{G}_j$, $\overline{G}_j$: Minimum and maximum generation at bus $j$
    \item[] $\gamma$: Index of the slack bus
    \item[] $L$: Economic cost  load loss per unit power per hurricane
    \item [] $\mathcal{H}_i^{max}, \mathcal{T}_i^{max}$: Largest element in set $\mathcal{H}_i$ and $\mathcal{T}_i$
    \item [] $R$: Time to restore power in the hurricane's aftermath
    \item[] $\prob(k)$: Probability that $\xi_1 = k$
    \item[] $\prob(k,m)$: Probability that $\xi_1 = k$ and $\xi_2 = m$
    \item [] $\Bar{\theta}_r$: Maximum difference in the phase angles of the endpoint buses of branch $r$
\end{itemize}
\noindent\textbf{Variables}
\begin{itemize}
    \item[] \textbf{First-stage variables}
    \item[] $n$: Integer variable, representing the number of Tiger Dam units that are acquired (each unit of Tiger Dam provides an additional level of resilience)
    \item[] $h_{il}$: Binary variable, indicating whether substation $i$ is hardened to level $l$
    \item[] \textbf{Second-stage variables}
    \item[] $t_{il}^k$: Binary variable, indicating whether Tiger Dam units are deployed at substation $i$ to level $l$ when $\xi_1 = k$
    \item[] \textbf{Third-stage variables}
    \item[] $z_{j}$: Binary variable, indicating whether bus $j$ is operational
    \item[] $u_{j}$: Binary variable, indicating whether generator at bus $j$ is dispatched
    \item[] $s_{j}$: Non-negative real variable, load satisfied at bus $j$
    \item[] $g_{j}$: Non-negative real variable, power generated at bus $j$
    \item[] $\theta_{j}$: Real variable, voltage phase angle of bus $j$
    \item[] $e_{r}$: Real variable, power flowing in branch $r$ from ${r.\mathrm{h}}$ to ${r.\mathrm{t}}$ where a negative value implies that the flow is in the opposite direction
\end{itemize}

In the above notation, we observe that the third stage variables $z_j, u_j, s_j, g_j, \theta_j, \text{and } e_r$ are indexed only in $j$ ($r$ in case of $e_r$). We note that these variables are defined within the recourse function and thus exist for all valid combinations of $k$ and $m$.

\subsection{Mathematical Formulation}\label{formulation}
In the rest of this article, we refer to the proposed three-stage stochastic MIP model as the Budget Allocation Model (BAM) with the mathematical formulation that follows: 
\begin{multline}
    \text{minimize } \sum_{i \in \mathcal{I}_f} \sum_{l \in \mathcal{H}_i}V_{il}h_{il}
       + An + Y\sum_{k \in \mathcal{P}} \prob(k) \sum_{i \in \mathcal{I}_f} \sum_{l \in \mathcal{T}_i}C_{il}t_{il}^k
       + \\ Y L \sum_{k \in \mathcal{P}} \sum_{m \in \mathcal{Q}_k} \prob(k,m) ~\mathcal{L}(\textit{\textbf{h}},\textit{\textbf{t}}^\textit{\textbf{k}},k,m),\label{bam_objective}
\end{multline}
\begin{align}
    & \text{subject to: \quad}  \sum_{l \in \mathcal{H}_i}h_{il} \leq 1, & \forall  i \in \mathcal{I}_f, \label{bam_hardening}\\
    & \sum_{l \in \mathcal{T}_i}t_{il}^k \leq  1, & \forall i \in \mathcal{I}_f, \quad k \in \mathcal{P},\label{bam_Tiger_dam}\\
    & \sum_{i \in \mathcal{I}_f} \sum_{l \in \mathcal{T}_i}lt_{il}^k \leq n, & \forall k \in \mathcal{P}.\label{bam_total_dams}
\end{align}

The objective function \eqref{bam_objective} in the BAM represents the total expected ``disaster management'' cost due to the flooding of power grid components during the planning horizon. The expression consists of four terms. The first term represents the cost of substation hardening. The second term refers to the acquisition cost of Tiger Dams. The third is the expected Tiger Dam deployment costs incurred during the planning horizon. The last term is representative of the economic costs associated with load loss and recovery operations. An explanation on how this value is computed is deferred for discussion until later in Section \ref{economics_discussion}.

Constraints \eqref{bam_hardening} are associated with the first-stage substation hardening decisions which ensure that each substation is assigned to exactly one hardening level. Similarly, constraints \eqref{bam_Tiger_dam}, which represent the second-stage Tiger Dam deployment decisions, ensure that each substation is assigned to exactly to one level of Tiger Dam deployment in the preparedness phase before each storm. Constraints \eqref{bam_total_dams} are coupling constraints which ensure that the total number of the Tiger Dam units deployed within any preparedness model in the second stage do not exceed the number of units acquired in the first stage.  

\textbf{Flood-aware DC power flow model:} The recourse function $\mathcal{L}(\textit{\textbf{h}},\textit{\textbf{t}}^\textit{\textbf{k}},k,m)$ in \eqref{bam_objective} is given by:
\begin{align}
\mathcal{L}(\textit{\textbf{h}},\textit{\textbf{t}}^\textit{\textbf{k}},k,m) = \text{minimize}  \sum_{j \in \mathcal{J}} D_j - s_{j}, \label{recourse_objective}
\end{align}
\begin{multline}
    \text{subject to:} \quad  \Delta_{ikm} \leq \max \left(\sum_{l \in \mathcal{H}_{i}}lh_{il}, \sum_{l \in \mathcal{T}_{i}}lt_{il}^k\right) + M(1-z_j),  \hfill \forall j \in \mathcal{B}_i, \quad \forall i \in \mathcal{I}_f,\label{bam_linking_1}
\end{multline}
\begin{multline}
\Delta_{ikm} \geq \max \left(\sum_{l \in \mathcal{H}_{i}}lh_{il}, \sum_{l \in \mathcal{T}_{i}}lt_{il}^k\right) -Mz_j + 1, \hfill \forall j \in \mathcal{B}_i, \quad \forall i \in \mathcal{I}_f,\label{bam_linking_2}
\end{multline}

\begin{align}
    & z_j = 1, & \forall j \in \mathcal{B}_i, \quad \forall i \not \in \mathcal{I}_f,\label{bam_linking_3}\\
    & u_{j} \leq z_{j}, &  \forall j \in \mathcal{J},\label{bam_dispatch}\\
    & s_{j} \leq z_{j}D_j, & \forall j \in \mathcal{J},\label{bam_supply_ub}\\
    & u_{j}\underline{G}_j \leq g_{j} \leq u_{j}\overline{G}_j, & \forall j \in \mathcal{J},\label{bam_generation_constraints}\\
    & -z_{r.\mathrm{h}} P_r \leq e_{r} \leq z_{r.\mathrm{h}} P_r, & \forall r \in \mathcal{R},\label{bam_edge_capacity_1}\\
    & -z_{r.\mathrm{t}} P_r \leq e_{r} \leq z_{r.\mathrm{t}} P_r, & \forall r \in \mathcal{R},\label{bam_edge_capacity_2}\\
    &|e_r - B_r (\theta_{r.\mathrm{h}} - \theta_{r.\mathrm{t}})|
                   \leq M(1 - z_{r.\mathrm{h}}z_{r.\mathrm{t}}), & \forall r \in \mathcal{R},\label{bam_phase_angle}\\
& \sum_{r \in \mathcal{N}_j^{out}}e_{r} - \sum_{r \in \mathcal{N}_j^{in}}e_{r}= g_{j} - s_{j}, & \forall j \in \mathcal{J}, \label{bam_flow_conservation}\\
& -\pi \leq \theta_{j} \leq \pi, &  \forall j \in \mathcal{J}, \label{bam_phase_side}\\
& \theta_{\gamma} = 0, \label{bam_slack_bus}\\
& |\theta_{r.\mathrm{h}} - \theta_{r.\mathrm{t}}| \leq \Bar{\theta}_r, & \forall r \in \mathcal{R} \label{bam_phase_bound},
\end{align}
where the first two inputs in the recourse function $\mathcal{L}(\textit{\textbf{h}},\textit{\textbf{t}}^\textit{\textbf{k}},k,m)$ represent the substation hardening and Tiger Dam deployment decisions, respectively. Specifically, $\textbf{\textit{h}} = (h_{il})_{i \in \mathcal{I}_f, l \in \mathcal{H}_i}$ and $\textbf{\textit{t}}^{\textbf{k}} = (t_{il}^k)_{i \in \mathcal{I}_f, l \in \mathcal{T}_i}$. 

\textbf{Uncertainty representation:} The last two elements in the input ($k$ and $m$) correspond to the value of the realizations of random variables $\xi_1$ and $\xi_2$, respectively. Each combination of the realized values represents a unique sample path corresponding to a flooding scenario with certain hurricane characteristics. Several such paths combine to form a scenario tree. We note that the flood scenarios associated with different sample paths can be identical. This is further discussed with an example in  Section \ref{section_scenario_tree}. In the rest of the article, we refer to scenarios using the tuple $(k,m)$ as their unique identifier. 

\textbf{Economic cost of load shedding:} The output of the recourse function $\mathcal{L}(\textit{\textbf{h}},\textit{\textbf{t}}^\textit{\textbf{k}},k,m)$ represents the total load loss in scenario $(k,m)$ when the substation hardening and Tiger Dam deployment decisions are \textit{\textbf{h}} and $\text{\textit{\textbf{t}}}^\text{\textbf{k}}$. An expectation over the output of the recourse function is taken for possible realizations of $\xi_1$ and $\xi_2$ in \eqref{bam_objective} to compute the expected load loss. The expected load loss value is then multiplied with appropriate constants to represent the economic costs associated with load loss and recovery operations. To understand how this is done, we first determine the value of the parameter $L$, which itself is the product of the economic value of load loss per hour (expressed in \$ per MW-hr), the restoration/down time which we assume that all the substations that are not operational experience, and a constant term representing the base-load value. The product $L \sum_{k \in \mathcal{P}} \sum_{m \in \mathcal{Q}_k} \prob(k,m) ~\mathcal{L}(\textit{\textbf{h}},\textit{\textbf{t}}^\textit{\textbf{k}},k,m)$ therefore represents the expected economic loss per hurricane. We finally multiply this value by $Y$ to compute the expected total economic and social costs incurred due to load losses in all the storms during the planning horizon. Further discussion on computing the economic cost of load loss for the state of Texas is presented in Section \ref{economics_discussion}.

\textbf{Contingency impact constraints:} In the definition of the recourse function $\mathcal{L}(\textit{\textbf{h}},\textit{\textbf{t}}^\textit{\textbf{k}},k,m)$, constraints \eqref{bam_linking_1} and \eqref{bam_linking_2} link the decisions made across multiple stages. Specifically, they compare the substation hardening level and the level of the Tiger Dam with the flood height at a substation in a particular scenario to determine if the substation is flooded. Accordingly, the value of $z_j$ is determined to be 0 or 1  where 0 indicates that the substation's protection is not sufficient and in turn the bus is not operational. For the substations that are not flooded in any scenario, the status of all the buses that they contain is set to being operational in constraints \eqref{bam_linking_3}. 

\textbf{Power flow constraints:} Constraints \eqref{bam_dispatch} capture the generator dispatch decisions. The constraints imply that a generator bus can only be dispatched when operational. Constraints \eqref{bam_supply_ub} enforce that if a bus is operational, the power supplied at the bus is at most the bus demand, and 0 otherwise. Similarly, constraints \eqref{bam_generation_constraints} ensure that when a generator is operational, the power dispatched from the generator is within its generation limits and is 0 otherwise.  Constraints \eqref{bam_edge_capacity_1} and \eqref{bam_edge_capacity_2} place limits on the power flow of branches and imply that if the bus at either end of the branch is not operational, then no power can flow through the branch. Constraints \eqref{bam_phase_angle} represent the Ohm's law constraint $B_r (\theta_{r.\mathrm{h}} - \theta_{r.\mathrm{t}}) =  e_{r}$ for each operational branch $r$. The big-$M$ based reformulation is used to impose the condition that the phase angle constraint need to only be enforced on a branch if both ends of the branch are operational. Constraints \eqref{bam_flow_conservation}  ensure flow conservation at each bus. Constraints \eqref{bam_phase_side} provide lower and upper limits of the phase angles, and constraint \eqref{bam_slack_bus} sets the phase angle for the reference or slack bus. Lastly, constraints \eqref{bam_phase_bound} limit the phase angle difference between endpoint buses of each branch.

We note that the mathematical formulation of the proposed three-stage model is a mixed integer non-linear optimization problem. This non-linearity is due to the presence of the $\max(\cdot)$ function in constraints \eqref{bam_linking_1} and \eqref{bam_linking_2}, the presence of  the absolute value function and bi-linear terms in constraints \eqref{bam_phase_angle}, and the presence of the absolute value function in constraints \eqref{bam_phase_bound}. The problem can, however, be reformulated into a mixed integer linear program using standard transformations. A detailed discussion on the reformulation of non-linear constraints is presented in Appendix \ref{reformulation}, and the determination of the values of big-$M$ is discussed in Appendix \ref{tightening_label_suppl}.



\subsection{Value of full coordination}\label{vofc}
While making the substation hardening decisions during the mitigation phase, the BAM accounts for the Tiger Dam-based preparedness measures that can be taken to protect the substations from flooding before an imminent hurricane. 
We now explain how we can quantify the benefit of decision making based on this integrated approach as opposed to making these decisions in a decoupled way.

Consider a case when the mitigation decisions are made without accounting for the preparedness options that might be available before an impending hurricane. To reflect this in the BAM, we set $n = 0$ and refer to this constrained version of the budget allocation model with no preparedness measures as the BAM-NP. 
Consequently, to decide the optimal budget for substation hardening, the BAM-NP will trade off hardening a substation vs leaving it unprotected from flooding which can subsequently cause load loss. Next, we fix the substation hardening levels to their respective values as obtained from the BAM-NP in the BAM. In this case, another three-stage version of the BAM will consider these substation hardening decisions as given and trade off the Tiger Dam acquisition and deployment costs with load loss to decide an optimal allocation of the budget for the purchase of Tiger Dams. We refer to this decoupled variant of the BAM as the BAM-D. We call the difference in the objective function values of the BAM-D and the BAM as the value of full coordination.

\section{Case study} \label{case_study}
This section presents a case study focused on the Texas coastal region to demonstrate the use of the BAM for addressing various questions related to the resilience of the power grid to extreme flood events. The input to the BAM consists of two main components. The first is a set of flood scenarios representing the impact of the different types of hurricanes hitting the Texas coast and the second is the
topology and the network parameters of the power grid that resembles the actual grid of the Texas coastal region. 
Further note that, we solve the BAM instances described in this case study on a high performance computing server with AMD's EPYC Milan processor (128 cores) and 256GB of RAM. The models are developed in gurobipy and solved using the Gurobi solver with the barrier algorithm \cite{gurobi}.

\subsection{Flood scenarios}
\begin{figure}
    \centering
    \includegraphics[scale=0.45]{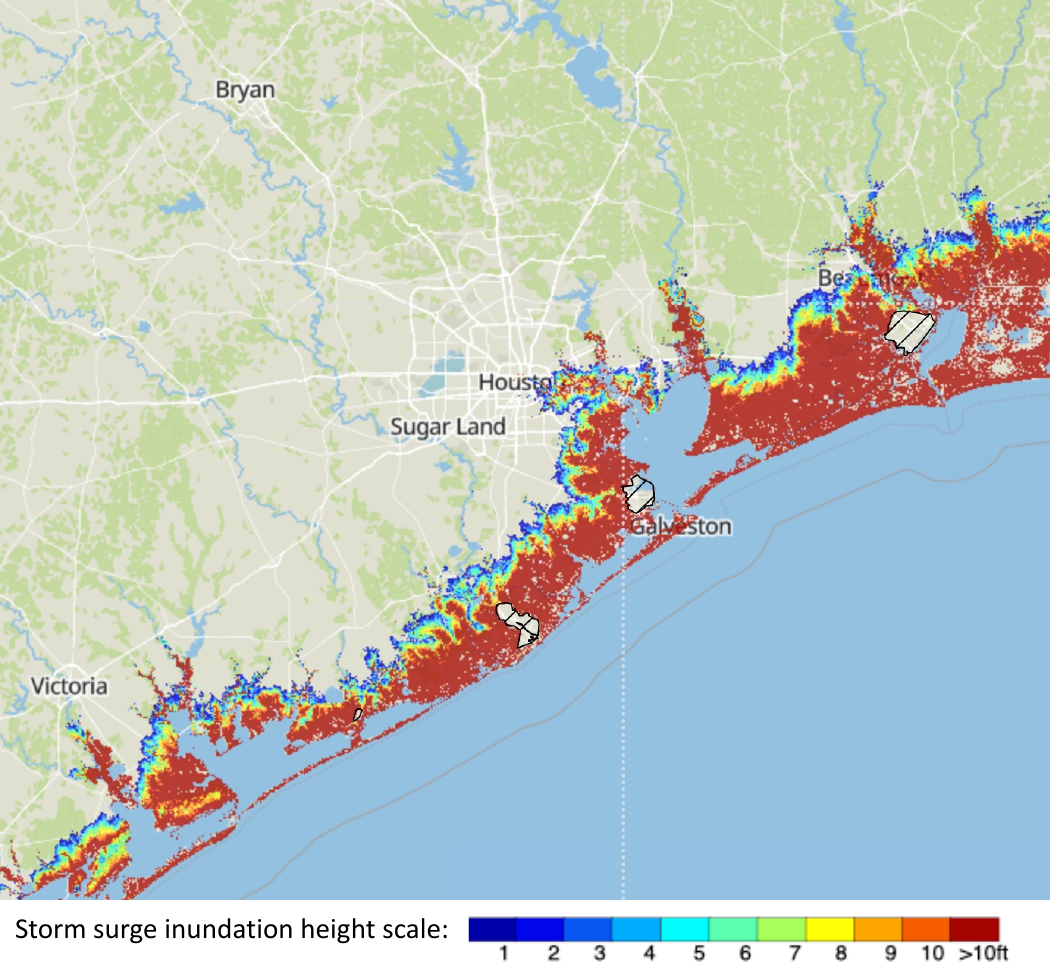}
    \caption{A sample MEOW generated using category 5 storms approaching the Texas Gulf Coast in the north-west direction with a forward speed of 5 mph}
    \label{fig:MEOW}
\end{figure}

We use the storm-surge maps \cite{MEOW} developed by the National Oceanic and Atmospheric Administration (NOAA) using the Sea, Lake, and Overland Surges from Hurricanes (SLOSH) model \cite{SLOSH} to sample the flood scenarios. The maps are generated by simulating hurricanes with different characteristics, such as intensity, forward speed, direction, and tide levels. For each combination of such parameters, multiple SLOSH runs (each with a different landfall location) are conducted to capture the variability in the storm surge. The maximum water levels from these runs are used to construct Maximum Envelopes of Water (MEOW) maps for different types of hurricanes. An example map for the Texas coastal region is shown in Figure \ref{fig:MEOW}. We overlay the Texas's synthetic coastal transmission grid on such a map to identify substations that are flooded due to inadequate protection. Accordingly, the buses within the substations and branches connected to those buses are considered out of order. This is modeled using the flood-aware dc power model with network reconfiguration constraints as described in Section \ref{formulation}. 

\subsection{Transmission grid}

\begin{figure}
    \centering
    \includegraphics[trim={0cm 0cm 0cm 0.9cm},clip, width=0.7\textwidth]{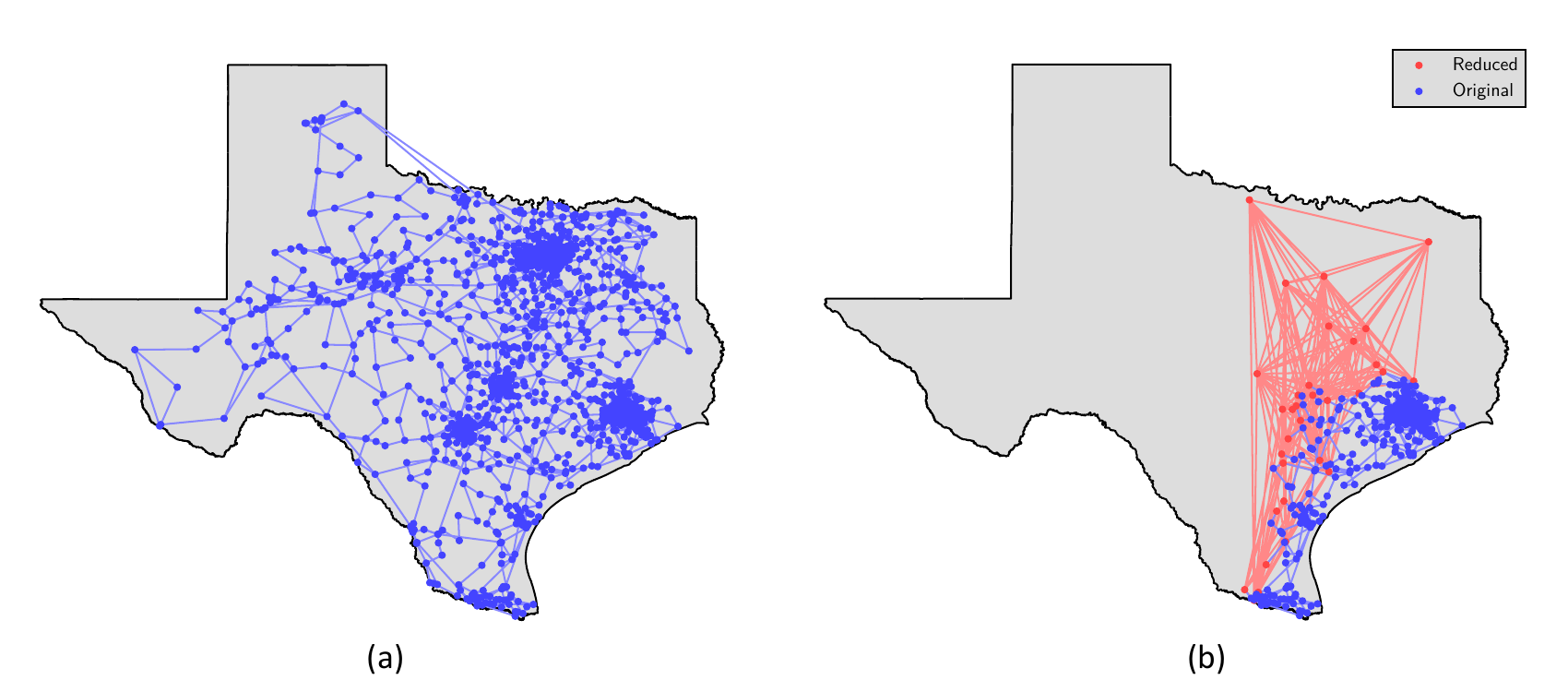}
    \caption{(a) The ACTIVSg2000 synthetic grid for Texas. (b) The reduced grid obtained after performing the network reduction. The red elements represent the new nodes and branches that were introduced as the artifacts of the reduction procedure to maintain equivalence in the grid characteristics}
    \label{fig:grid_reduction}
\end{figure}

We use the synthetic transmission grid ACTIVSg2000 developed as part of an ARPA-E project \cite{activsg2000} on the footprint of the actual Texas grid. The synthetic grid consists of 2000 buses distributed across 1250 substations and 3206 branches. Moreover, to keep the case study realistic while maintaining computational tractability, we make two modifications to ACTIVSg2000 dataset to assess and improve its resilience to coastal storm surge risk. First, the coordinates of the 1250 substations are remapped using real-world substation locations obtained from the Homeland Infrastructure Foundation-Level Data (HIFLD) Electric Substations dataset \cite{HIFLD_Electric_Substations}. An optimization problem is solved to minimize the displacement caused by the remapping, ensuring that the electrical structure of the power grid remains unchanged. These modifications allow for a more realistic representation of the grid while considering real-world flood risks via NOAA-based flood scenarios. Second, a network reduction is performed using the electrical equivalent (EEQV) feature in {PSS\textregistered E} to retain only the grid components susceptible to storm surge-induced flooding. Buses that are in the inland region and not exposed to storm surge are aggregated into a smaller set of nodes, while the coastal part of the grid remains almost unchanged. The topological changes resulting from the reduction are illustrated in Figure \ref{fig:grid_reduction}. The obtained reduced grid consists of all of the 362 substations from the coastal grid, encompassing 663 buses and 1509 branches, with a generation capacity of 60.98 GW and a total load of 39.69 GW.

\subsection{Scenario tree} \label{section_scenario_tree}
The NOAA storm-surge maps are developed using flood simulations while considering eight different storm directions (west-south-west, west, west-north-west, north-west, north-north-west, north, north-north-east, and north-east), six different intensity categories (0-5), and four different forward speeds (5, 10, 15, and 25 mph). Consequently, the dataset comprises of 192 unique flood maps, one for each combination of storm direction, category, and forward speed.

\begin{figure}
    \centering
    \includegraphics[width=0.6\textwidth]{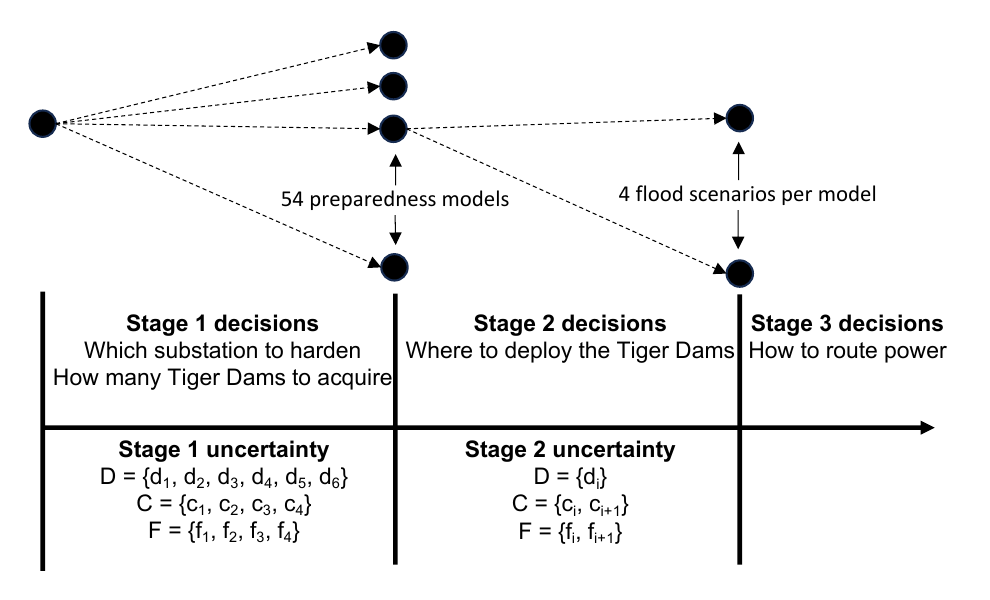}
    \caption{A schematic diagram representing the scenario tree used in the case study. Here, D, C, and F represent the set of storm directions, intensity categories, and forward speeds.}
    \label{fig:scenario_tree}
\end{figure}

To demonstrate the usefulness of the BAM with a computationally tractable use case, we reduce the size of the problem by eliminating a subset of less severe scenarios. Specifically, we first drop the scenarios corresponding to two directions (west-south-west and northeast) as hurricanes belonging to these categories do not cause significant flooding in the Texas Gulf Coast. We also drop the hurricanes of categories 0 and 1. The hurricanes belonging to these categories cause little to no flooding at the substations considered in this study. As a result, we have 96 flood scenarios belonging to 6 directions and four different levels of hurricane categories and forward speeds. These 96 scenarios are considered likely to impact the power grid and are part of the uncertainty set at the ``root node'' of the scenario tree we consider (Figure \ref{fig:scenario_tree}). At the root node representing the first stage of decision making, we make the here-and-now substation hardening and Tiger Dam acquisition decisions that are to mitigate against the flooding risks in all 96 hurricane scenarios.

Before we make the second-stage Tiger Dam deployment decisions, the direction of the hurricane is revealed and the uncertainty about the category and the forward speed of the hurricane is reduced. This is represented in Figure \ref{fig:scenario_tree} where both the category and forward speed can take one of the two consecutive values from the ordered set representing the forward speed and category. Consequently, in the second stage, we consider 54 preparedness models (for each combination $d_i, c_i, c_{i+1}, f_i, f_{i+1}$), each consisting of 4 scenarios (for each combination of $d_i, c_i, f_i$) while making the Tiger Dam deployment decisions. Subsequently, a scenario is realized causing a subset of substations with inadequate protection to go out of order. In this study, we assume the path probabilities to be equal in the case study. However, the model is flexible, in that a decision maker can use their own probabilistic beliefs.

In the third stage, we make decisions about re-routing power with the grid components remaining operational such that the load loss is minimized. 

\subsection{Parametric Setup}\label{setup}
\subsubsection{Substation hardening} For each substation $i$ that is flooded in at least one scenario, we define $\mathcal{H}_i = \{3, 6, \cdot \cdot \cdot, 21\}$. The value of $V_{il}$, which represents the cost of hardening substation $i$ to resilience level $l$ is considered to be directly proportional to $l$ and is given by \$100,000 times $l$. This implies that it will take at most \$2.1M to fully harden a substation against any flooding as the maximum flood height across all the scenarios is 21 feet. This cost is close to the approximated value of \$2.2M in a policy brief published by the Bass Connections at Duke University \cite{substation_cost}. The hardening cost, however, can vary depending on the location of the substation, terrain, wall thickness and other related specifications. Therefore to study the impact that a change in the substation hardening cost has on the optimal budget allocation, we do a sensitivity analysis by varying the value of $V_{il}$ from \$25,000/ft to \$100,000/ft at \$25,000 intervals in Section \ref{bam_variable_cost}.

\subsubsection{Tiger Dams} Similar to substation hardening, we define $\mathcal{T}_i = \{3, 6\}$ as the set representing the Tiger Dam deployment levels. To further assess the impact that the maximum Tiger Dam deployment level ($\mathcal{T}_i^{\text{max}}$) has on the optimal budget allocation, we vary its value from 6 feet in the default case to 21 feet at 3 feet intervals in Section \ref{bam_tiger_dam_improvement}.

The cost of each Tiger Dam unit is assumed to be \$40,000/ft. Therefore, Tiger Dam units worth \$120,000 and \$240,000 are required to build ad hoc protection barriers of height 3 and 6 feet, respectively. These values are based on the assumption that a Tiger Dam is expected to surround the periphery of a substation which is considered to be a square with each side measuring 200 feet. The Tiger Dams cost around \$200 per 4-foot piece. Therefore, to protect the perimeter of 800 feet, we need 200 pieces which will cost \$40,000. Lastly, the cost incurred for each Tiger Dam deployment is assumed to be \$10,000. 

\subsubsection{Economic and social cost of load loss}\label{economics_discussion}
The expected cost due to load loss in a hurricane is determined using three values: the value of load loss (VOLL), the expected time to restore power ($R$) , and the expected load loss per hurricane. The parameter $L$ in the objective function represents the product of the first two terms times 100; the constant 100 is to convert the output load loss from the recourse function $\mathcal{L}(\textit{\textbf{h}},\textit{\textbf{t}}^\textit{\textbf{k}},k,m)$, which is in per unit system, to megawatt-hours. The default value of VOLL is assumed to be \$6,000 per MW-hr and is in line with the value approximated during a study conducted by the Texas grid operator, the Electric Reliability Council of Texas (ERCOT) \cite{ercot_voll}. The default value of $R$ is assumed to be 48 hours. Therefore, the value of $L$ that we use in our study as the base case is \$2.88M per unit of power per hurricane. Moreover, the expected number of hurricanes during the planning horizon is assumed to be 10. Therefore, the last term in the objective function is reduced to \$28.8M times $\sum_{k \in \mathcal{P}} \sum_{m \in \mathcal{Q}_k} \prob(k,m) ~\mathcal{L}(\textit{\textbf{h}},\textit{\textbf{t}}^\textit{\textbf{k}},k,m)$ with the default values of the parameters. To further assess the impact of the changes in VOLL and power restoration times, a parametric study is conducted in Section \ref{voll_sensitivity} and Section \ref{time_sensitivity}.

\section{Results and Discussion}
\subsection{A note on figures}\label{figure_discussion}
Figures \ref{fig:granular_voll_250}, \ref{fig:granular_voll_6000}, \ref{fig:decoupled_6000}, and \ref{fig:max_prep_15} are all color coded rectangular grids where each row represents the permanent hardening and Tiger Dam deployment decisions at a substation. In all the figures, we consider the same set of substations placed in the same sequence. They are sorted from top to bottom in the increasing order of average flood values. Moreover, each substation is flooded in at least one of the 96 distinct flooding scenarios considered in this study. The left-most column in the figures is marked with only one of the two colors (black and yellow) and represents the first-stage permanent hardening decisions. The next 54 columns show the Tiger Dam deployment decisions in each of the 54 preparedness models associated with different realizations of $\xi_1$ as  explained in Section \ref{section_scenario_tree}.  

A grid cell with index $(i,k)$ representing substation $i$ and preparedness model $k$ in the aforementioned rectangular grid is colored white when the variable max\_flood, which denotes the maximum flooding level across all preparedness phase scenarios and $t_{il}^k, ~ \forall l \in \mathcal{T}_i$ are all simultaneously 0. That is, no flooding occurs at substation $i$ in any of the scenarios considered within the preparedness model associated with $\xi_1 = k$ and therefore, no Tiger Dams are deployed. 

When the cell is blue, it signifies that max\_flood $>$ 0 and $t_{il}^k = 0, \forall l \in \mathcal{T}_i$. It implies that although the substation is flooded in at least one scenario, Tiger Dams are not deployed at the substation. This can be due to two reasons. The first is that the substation has been chosen for permanent hardening in the first stage (as is the case in most of the substations). The second is that because we have limited Tiger Dam resources, it is deemed more cost-effective to leave the substation unprotected and subject it to the risk of flooding while using the limited Tiger Dam resources to protect other substations which are more likely to contribute to minimizing the total load loss induced cost. If the cell is green, a Tiger Dam is deployed in the preparedness phase and is capable of safeguarding the substation against the maximum flooding (which is non-negative) that can occur in one or more of the scenarios considered within a particular preparedness model. Lastly, orange color is indicative of substations where a Tiger Dam provides protection in some scenarios but not all as the maximum flood level exceeds its protective capacity. 

\subsection{Sensitivity to VOLL} \label{voll_sensitivity}
To assess how the allocation of the optimal budget changes with VOLL, we solve the BAM for different VOLL as shown in Table \ref{tab:coupled_voll}. In the table, we observe that the total disaster cost increases with the increase in VOLL. This can be attributed to the increase in the budget allocated for hardening of the substations from \$33M to \$64.2M as VOLL increases from \$250 per MW-hr to \$6,000 per MW-hr. When VOLL is \$250 per MW-hr, the BAM protects a smaller set of substations using permanent hardening. In this case, since the penalty on load loss is relatively small, the BAM recommends protecting several substations using Tiger Dams when flooding is lower than the maximum flood protection that Tiger Dams can provide and leave them unprotected to higher flood levels. 

\begin{figure}
  \centering

  \begin{subfigure}{0.45\textwidth}
    \centering
    \includegraphics[trim={0cm 0cm 4.9cm 0cm},clip, height=2.5in]{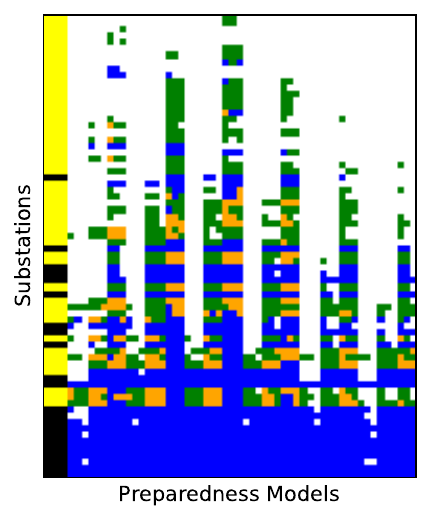}
    \caption{\$250/MWh}
    \label{fig:granular_voll_250}
  \end{subfigure}
  \begin{subfigure}{0.45\textwidth}
    \centering
    \includegraphics[trim={0.55cm 0cm 4.9cm 0cm},clip, height=2.5in]{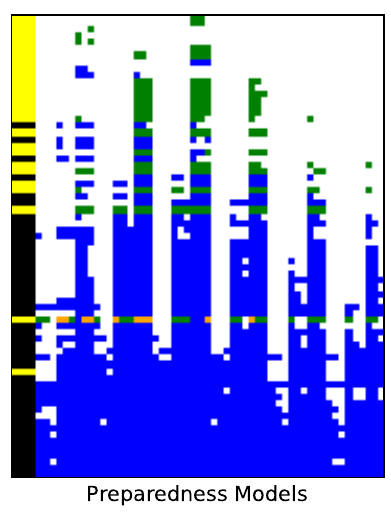}
    \caption{\$6,000/MWh}
    \label{fig:granular_voll_6000}
  \end{subfigure}

  \vspace{\baselineskip}

  \begin{subfigure}{0.45\textwidth}
    \centering
    \includegraphics[trim={0cm 0cm 4.9cm 0cm},clip, height=2.5in]{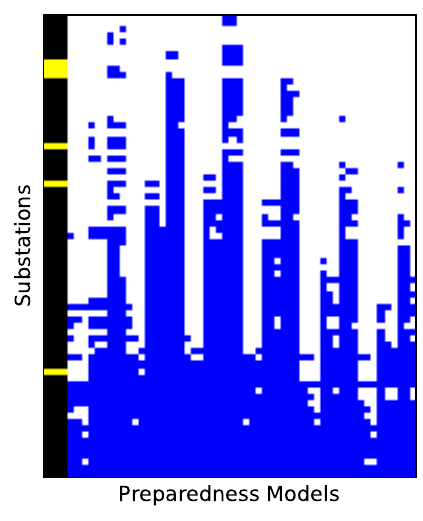}
    \caption{\$6000/MWh (BAM-D)}
    \label{fig:decoupled_6000}
  \end{subfigure}
  \begin{subfigure}{0.45\textwidth}
    \centering
    \includegraphics[trim={0.55cm 0cm 0.0cm 0cm},clip, height=2.5in]{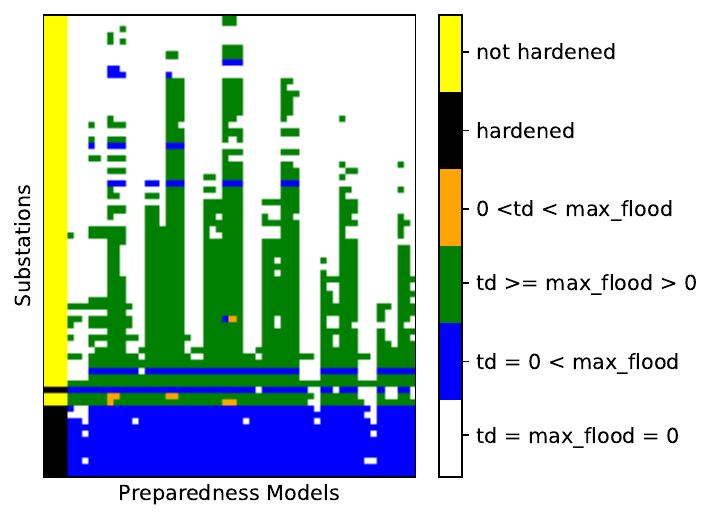}
    \caption{15 feet}
    \label{fig:max_prep_15}
  \end{subfigure}

  \caption{First- and second-stage solutions: \eqref{fig:granular_voll_250} and \eqref{fig:granular_voll_6000} depict the solutions of the BAM when VOLL is \$250/MW-hr and \$6000/MW-hr, respectively. \eqref{fig:decoupled_6000} shows the same for the BAM-D for a VOLL of \$6000/MW-hr. \eqref{fig:max_prep_15} represents the BAM solutions for the maximum preparedness level of 15 feet. For more details, refer to Section \ref{figure_discussion}}
  \label{fig:overall_label}
\end{figure}

A higher VOLL implies a larger penalty on the loss of load, and therefore the model minimizes that component of the total cost by hardening increasingly more substations to eventually attain near-zero load loss at a VOLL of \$6,000 per MW-hr. To verify this, refer to Figures \ref{fig:granular_voll_250} and \ref{fig:granular_voll_6000}. We further notice another trend: increasing the hardening budget generally leads to a reduction in the budget allocated for acquisition and deployment of Tiger Dams. This is expected as permanent hardening of additional substations obviates the need for Tiger Dam that are otherwise deployed to protect them. Consequently, this further leads to a reduction in the expected deployment costs.

\begin{table}
    \centering
    \caption{Budget allocation for substation hardening and Tiger Dam acquisition for the different values of load loss, VOLL. Also shown are the expected costs of Tiger Dam deployment and of load loss in the planning horizon.}
    \begin{tabular}{|c|c|c|c|c|c|}
    \hline
    VOLL &  Hardening &  Tiger Dam &  Deployment &  Load Loss &  Total \\
    \hline
    250  &       33.0 &       4.92 &        1.35 &      11.32 &  50.59 \\
    500  &       45.0 &       3.60 &        0.94 &       7.87 &  57.41 \\
    1000 &       51.3 &       3.60 &        0.67 &       5.48 &  61.05 \\
    2000 &       56.4 &       3.48 &        0.51 &       5.20 &  65.59 \\
    3000 &       59.7 &       3.48 &        0.42 &       3.50 &  67.10 \\
    4000 &       62.1 &       3.24 &        0.39 &       2.32 &  68.05 \\
    5000 &       62.4 &       3.36 &        0.39 &       2.42 &  68.57 \\
    6000 &       64.2 &       3.12 &        0.33 &       1.37 &  69.02 \\
    \hline
    \end{tabular}
    \label{tab:coupled_voll}
\end{table}



In the last row of Table \ref{tab:coupled_voll}, we observe that the optimal budget for substation hardening is \$64.2M. This is for a VOLL of \$6,000 per MW-hr, a value close to what has been approximated for the state of Texas in a study by the grid operator, the Electric Reliability Council of Texas (ERCOT) \cite{ercot_voll}. Therefore, in this case, we should prioritize making significant investments in permanent hardening of the substations located in the Texas coastal grid to protect them against storm surge induced flooding. The other point to emphasize is that for any VOLL, there exists an implicit optimal budget for substation hardening and Tiger Dam acquisition that minimizes the expected overall disaster cost. If the allocated mitigation budget is any smaller than this value, it will only increase the overall disaster cost.

\subsection{Sensitivity to restoration time} \label{time_sensitivity}
In Section \ref{economics_discussion}, we have discussed that the expected total load loss induced cost is given by $Y L \sum_{k \in \mathcal{P}} \sum_{m \in \mathcal{Q}_k} \prob(k,m) \mathcal{L}(\textit{\textbf{h}},\textit{\textbf{t}}^\textit{\textbf{k}},k,m)$ where the parameter $L$ is a constant times the product of VOLL and $R$. Further recall that the default value of VOLL and $R$ is \$6000 per MW-hr and 48 hours, respectively. Therefore, for a VOLL of  \$2000 per MW-hr, the effect of change in the restoration time from 48 to 24 hours is reflected in the solution corresponding to the VOLL of \$1000 per MW-hr in Table \ref{tab:coupled_voll}. Similarly, to obtain the budget allocation when $R$ is reduced to 12 or 6 hours, we refer to the budget associated with VOLL of \$500 and \$250 per MW-hr in Table \ref{tab:coupled_voll}, respectively. From the table, we further infer that a decrease in restoration time will reduce both expected hardening and total disaster cost as power can be restored quickly and cause relatively less economic damage. It is also interesting to note that for the VOLL of \$6000 per MW-hr, which is close to what has been approximated for the state of Texas in \cite{ercot_voll}, the proportion of the total budget that is allocated to Tiger Dam acquisition and deployments remains almost unaffected even if the restoration time is reduced from the default value of 48 hours to 4 hours (this solution corresponds to the second row in Table \ref{tab:coupled_voll}). 

\subsection{Value of coordinated decision-making}
In Section \ref{vofc}, we introduce the concept of the value of coordination to quantify the benefit of accounting for Tiger Dam acquisition and deployment decisions while allocating the budget for substation hardening. We now compare the total disaster costs obtained from the BAM and the BAM-D in Table \ref{tab:coupled_voll} and \ref{tab:decoupled_voll}, respectively, for different values of load loss. We observe that the total disaster costs in the case of decoupled decision making performed solving the BAM-NP followed by the BAM-D is considerably more (between 9.4\% to 14.5\% depending on VOLL) than the costs of coupled decision making performed by solving the original (coordinated) BAM. 

In the decoupled case, because the substation hardening decisions do not account for preparedness measures, they tend to harden the substations which could have been protected using Tiger Dams. Due to this over-hardening of substations, the cost of Tiger Dam acquisition and deployment is much lower in the decoupled case as compared to the coupled case. This can also be confirmed by comparing the substation hardening decisions between Figure \ref{fig:granular_voll_6000} and Figure \ref{fig:decoupled_6000}, both for a VOLL of \$6000 per MW-hr. The plots for other VOLLs are in Appendix \ref{visualization_supplement}. Here, we emphasize that by not coordinating the decisions made in the hardening and preparedness phases, we end up with significant overspending even when the decisions in individual phases are optimized. When the VOLL is \$6,000 per MW-hr, the overspending in hardening is more than \$10M, although the net expected overspending is \$6.3M, achieved due to lower costs in Tiger Dams and load loss, both possible thanks to over hardening.


\begin{table}
    \centering
    \caption{Budget allocation for substation hardening and Tiger Dam acquisition for the different values of load loss, VOLL for the case of decoupled decision-making.}
    \begin{tabular}{|c|c|c|c|c|c|}
    \hline
    VOLL &  Hardening &  Tiger Dam &  Deployment &  Load Loss &  Total \\
    \hline
    250  &       47.1 &       1.44 &        0.53 &       8.86 &  57.93 \\
    500  &       59.4 &       0.84 &        0.22 &       4.65 &  65.11 \\
    1000 &       64.8 &       0.48 &        0.08 &       4.44 &  69.80 \\
    2000 &       69.9 &       0.36 &        0.07 &       2.48 &  72.81 \\
    3000 &       72.3 &       0.24 &        0.02 &       1.17 &  73.73 \\
    4000 &       73.5 &       0.24 &        0.01 &       0.92 &  74.67 \\
    5000 &       74.4 &       0.12 &        0.00 &       0.52 &  75.04 \\
    6000 &       74.7 &       0.00 &        0.00 &       0.62 &  75.32 \\
    \hline
    \end{tabular}
    \label{tab:decoupled_voll}
\end{table}

\subsection{Effect of improvement in short-term measures} \label{bam_tiger_dam_improvement}
We now assess the impact of an improvement in the efficacy of Tiger Dams on the optimal budget allocation. To do so, we vary the value of the maximum flood resilience that Tiger Dams can provide from 6 feet (which we have assumed so far) to 21 feet at 3 feet intervals. The optimal budget allocation for these different resilience levels is reported in Table \ref{tab:preparedness_phase}. We observe that as the effectiveness of preparedness measures improve, both the total disaster costs and the proportion allocated to permanent hardening decrease significantly. Specifically, comparing the last rows in Table \ref{tab:coupled_voll} and \ref{tab:preparedness_phase}, we see that if Tiger Dams are as effective as a permanent hardening measure like walls, it obviates the need for substation hardening and reduces the total cost from \$69.02M to \$29.33M. While this reduction is unrealistic unless there is a technological breakthrough in ad hoc barrier design and manufacturing technology, it serves as a bound on how low the cost could be if they were that effective. To visualize the differences in the solutions corresponding to first- and second-stage decisions between the base case when $\mathcal{T}_i^{\text{max}} = 6$ and the case when  $\mathcal{T}_i^{\text{max}} = 15$, we refer to Figures \ref{fig:granular_voll_6000} and \ref{fig:max_prep_15}. From the figures, we can infer that because Tiger Dams can effectively prevent flooding up to 15 feet, much fewer substations undergo permanent hardening. The solutions for the cases when $\mathcal{T}_i^{\text{max}}$ is set to 9,12,18, and 21 feet, respectively are in Appendix \ref{visualization_supplement}.
\begin{table}
    \centering
    \caption{Budget allocation for substation hardening and Tiger Dam acquisition for different values of maximum flood height that Tiger Dam can prevent.}
    \begin{tabular}{|c|c|c|c|c|c|}
    \hline
    Height &  Hardening &  Tiger Dam &  Deployment &  Load Loss &  Total \\
    \hline
    6    &       64.2 &       3.12 &        0.33 &       1.37 &  69.02 \\
    9    &       52.8 &       6.60 &        0.75 &       2.45 &  62.60 \\
    12   &       34.8 &      13.08 &        1.60 &       2.25 &  51.73 \\
    15   &       24.0 &      17.04 &        2.08 &       1.33 &  44.45 \\
    18   &       16.8 &      19.80 &        2.46 &       0.32 &  39.38 \\
    21   &        0.0 &      26.40 &        3.23 &       0.30 &  29.93 \\
    \hline
    \end{tabular}
    \label{tab:preparedness_phase}
\end{table}

\subsection{Effect of change in substation hardening cost}\label{bam_variable_cost}
In Table \ref{tab:variable_hardening}, we observe that as long as the per-foot variable hardening cost is more than the Tiger Dam unit cost (\$40,000), the recommended budget that should be allocated to Tiger Dam acquisition and deployment remains almost the same. We further notice that irrespective of how inexpensive permanent hardening is, it forms the major component of the overall disaster cost in an optimal budget. When the hardening budget is reduced to half of the default value of \$100,000/ft, the load loss cost is almost one fourth the original value, whereas if the hardening budget is reduced to 1/4 of its default value, the load loss cost decreases to approximately 1/27 of its original value. However, because the cost of load loss is much lower than the other costs (hardening and preparedness measures), its contribution to the change in the overall cost remains minimal as the hardening cost changes. Note that the cost of load loss is small as it is minimized via optimization because of the high load loss penalty at VOLL of \$6,000.  
 \begin{table}
    \centering
    \caption{Budget allocation as a function of the permanent hardening cost}
    \begin{tabular}{|c|c|c|c|c|c|}
    \hline
    V &  Hardening &  Tiger Dam &  Deployment &  Load Loss &  Total \\
    \hline
    25   &      18.45 &       0.24 &        0.02 &       0.05 &  18.76 \\
    50   &      32.85 &       3.12 &        0.27 &       0.35 &  36.59 \\
    75   &      48.38 &       3.12 &        0.33 &       1.10 &  52.93 \\
    100  &      64.20 &       3.12 &        0.33 &       1.37 &  69.02 \\
    \hline
    \end{tabular}
    \label{tab:variable_hardening}
\end{table}
\section{Conclusions}
We propose a three-stage stochastic MIP model to assist in transmission grid resilience planning against extreme flood events. We further develop a case study using a grid instance that represents the coastal part of the Texas's transmission grid and the storm-surge flood maps developed by NOAA. Using the flood maps as scenarios in the three-stage model, we capture the network-wide impact of correlated flooding on the grid. 

With the help of the case study, we demonstrate how the proposed model can be used to address a diverse range of questions related to grid resilience. This involves examining the influence of load loss value on the optimal budget for mitigation and preparedness measures. Additionally, we evaluate how factors such as the restoration time, the effectiveness of preparedness measures, and the cost of substation hardening impact the allocation of the optimal budget for power grid resilience. The significance of considering preparedness measures while making mitigation phase decisions is quantified by the model. The extensive parametric study reveals the necessity of substantial investments in substation hardening in addition to ad hoc flood barriers even when the loss of load is assumed to have a very low value. Failure to make such resilience investments can result in a disproportionate increase in economic losses in the aftermath of flood events. 

For future research, we propose developing models that account for equity while making budget allocation decisions for power grid resilience. Another research direction is the development of flood-scenario generation approaches that simultaneously account for both storm-surge and precipitation-induced flooding. 


\if0\blind{
\section*{Acknowledgements}
The authors acknowledge the generous support from The University of Texas at Austin's Energy Institute.	} \fi

\bibliographystyle{chicago}
\spacingset{1}
\bibliography{main}

\clearpage
\appendix
\numberwithin{figure}{section}
\section{Appendix}

\subsection{Reformulation of the non-linear constraints}\label{reformulation}

We first consider the linking constraints  \eqref{bam_linking_1} and \eqref{bam_linking_2}. These constraints can be linearized using the big-$M$-based reformulation. To do so, we observe that in an optimal solution, the Tiger Dam protection level at all of the flooded substations within each short-term model will either be 0 or a non-negative value that is greater than the permanent hardening level at that substation. This is because in all other cases, the deployment of a Tiger Dam will incur an additional deployment cost and yet will not provide any additional resilience as compared to what has been achieved with permanent hardening. 

Consequently, constraints \eqref{bam_linking_1} and \eqref{bam_linking_2} can be reformulated as follows. First, for simplicity, we replace  $\max \left(\sum_{l \in \mathcal{H}_{i}}lh_{il}, \sum_{l \in \mathcal{T}_{i}}lt_{il}^k\right)$ with $p_{ik}$ in both the constraints. Therefore, we can replace constraints \eqref{bam_linking_1} and \eqref{bam_linking_2} in the BAM with the following set of constraints:
\begin{subequations}
  \begin{align}
    \Delta_{ikm} \leq p_{ik} + M(1-z_j), \quad \forall j \in \mathcal{B}_i, \quad \forall i \in \mathcal{I}_f, \label{modified_linking_1}\\
    \Delta_{ikm} \geq p_{ik} -Mz_j + 1, \quad \forall j \in \mathcal{B}_i, \quad \forall i \in \mathcal{I}_f,\label{modified_linking_2}\\
    p_{ik} \geq \sum_{l \in \mathcal{H}_{i}}lh_{il}, \quad \forall i \in \mathcal{I}_f,\label{modified_linking_3}\\
    p_{ik} \geq \sum_{l \in \mathcal{T}_{i}}lt_{il}^k, \quad \forall i \in \mathcal{I}_f,\label{modified_linking_4}\\
    p_{ik} \leq \sum_{l \in \mathcal{H}_{i}}lh_{il} + M\sum_{l \in \mathcal{T}_i}t_{il}^k, \quad \forall i \in \mathcal{I}_f,\label{modified_linking_5}\\
    p_{ik} \leq \sum_{l \in \mathcal{T}_{i}}lt_{il}^k + M(1-\sum_{l \in \mathcal{T}_i}t_{il}^k),\quad \forall i \in \mathcal{I}_f\label{modified_linking_6}. 
  \end{align}
\end{subequations}
Moreover, based on the observation that if a Tiger Dam is deployed at a substation, its resilience level will always be higher than the permanent hardening level at that substation, we further employ the following set of valid inequalities: 
\begin{equation}
    h_{i,l'} + t_{i,l}^k \leq 1, \quad \forall l'\geq l, \quad i \in \mathcal{I}_f.
\end{equation}
Next, we linearize constraints \eqref{bam_phase_angle}. To do so, let us first introduce an additional variable $\alpha_r$ for each branch $r$ such that $\alpha_r = z_{r.\mathrm{h}}z_{r.\mathrm{t}}$. In this case, constraints \eqref{bam_phase_angle} is equivalent to the following:
\begin{subequations}
\begin{gather}
    |e_r - B_r (\theta_{r.\mathrm{h}} - \theta_{r.\mathrm{t}})| \leq M(1 - \alpha_r), \quad \forall r \in \mathcal{R},\label{bam_absloute_phase}\\
    \alpha_r \leq z_{r.\mathrm{h}}, \quad \forall r \in \mathcal{R},\\
    \alpha_r \leq z_{r.\mathrm{t}}, \quad \forall r \in \mathcal{R},\\
    \alpha_r \geq z_{r.\mathrm{h}} + z_{r.\mathrm{t}} - 1, \quad \forall r \in \mathcal{R}.           
\end{gather}
\end{subequations}
To linearize the absolute value function in constraints \eqref{bam_absloute_phase}, we replace it with the following constraints:
\begin{subequations}
    \begin{gather}
        e_r - B_r (\theta_{r.\mathrm{h}} - \theta_{r.\mathrm{t}})
                   \leq M(1 - \alpha_r), \quad \forall r \in \mathcal{R}\\
                   B_r (\theta_{r.\mathrm{h}} - \theta_{r.\mathrm{t}}) - e_r
                   \leq M(1 - \alpha_r), \quad \forall r \in \mathcal{R}.
    \end{gather}
\end{subequations}
Lastly, we replace constraints \eqref{bam_phase_bound} with the following linear equivalents:
\begin{subequations}
    \begin{gather}
\theta_{r.\mathrm{h}} - \theta_{r.\mathrm{t}} \leq \Bar{\theta}_r, \quad \forall r \in \mathcal{R}\\
\theta_{r.\mathrm{t}} - \theta_{r.\mathrm{h}} \leq \Bar{\theta}_r, \quad \forall r \in \mathcal{R}.
    \end{gather}
\end{subequations}

\subsection{Tightening the big-$M$ values in the BAM}\label{tightening_label_suppl}
The reformulated constraints in Section \ref{reformulation} have several constraints with big-$M$ terms which can be tightened further.
We claim that the smallest value of big-$M$ which satisfies constraints \eqref{modified_linking_1} is $W_i$. To verify this, observe that $-\max(\mathcal{H}_i^{max},\mathcal{T}_i^{max}) \leq \Delta_{ikm} - p_{ik} \leq W_{i}.$ Now, in constraints \eqref{modified_linking_1}, we need
$$\frac{\Delta_{ikm} - p_{ik}}{M} \leq 1. $$
The smallest value of $M$ that ensures this is $W_{i}$ and is for the case when $p_{ik} = 0$ and $\Delta_{ikm}$ takes the maximum value of $W_{i}$. Similarly, the big-$M$ value for constraints \eqref{modified_linking_2} is determined to be $\max(\mathcal{H}_i^{max},\mathcal{T}_i^{max}) + 1$. 

Next, we determine the value of big-$M$s for constraints  \eqref{modified_linking_5} and \eqref{modified_linking_6}. To do so, let us consider a case when $\sum_{l \in \mathcal{T}_i}t_{il}^k = 1$. Therefore, constraints \eqref{modified_linking_5} - \eqref{modified_linking_6} reduce to:
\begin{subequations}
  \begin{align}
    p_{ik} \leq \sum_{l \in \mathcal{H}_{i}}lh_{il} + M, \quad \forall i \in \mathcal{I}_f,\label{big_m_determination_3}\\
    p_{ik} \leq \sum_{l \in \mathcal{T}_{i}}lt_{il}^k,\quad \forall i \in \mathcal{I}_f\label{big_m_determination_4}. 
  \end{align}
\end{subequations}

In equation \eqref{big_m_determination_3}, a value of $M = \mathcal{T}_i^{max}$ ensures that constraints \eqref{big_m_determination_3} remain inconsequential and $p_{ik} = \sum_{l \in \mathcal{T}_{i}}lt_{il}^k$. Using the same approach, the value of $M$ in constraints \eqref{big_m_determination_4} is determined to be $\mathcal{H}_i^{max}$. 

Lastly, the value of $M$ in constraints \eqref{bam_absloute_phase} is $B_r\Bar{\theta}_r$. To verify this, consider the case when $\alpha_r = 0$. This happens when either one or both the end buses of branch $r$ are not operational. Consequently, constraints \eqref{bam_edge_capacity_1} and \eqref{bam_edge_capacity_2} imply that $e_r = 0$.  Therefore, we need to determine a value of $M$ that upper bounds the value of the expression $|B_r (\theta_{r.\mathrm{h}} - \theta_{r.\mathrm{t}})|$. We know that $| \theta_{r.\mathrm{h}} - \theta_{r.\mathrm{t}}| \leq \Bar{\theta}_r$ from constraints \eqref{bam_phase_bound} and hence, it follows that $M = B_r\Bar{\theta}_r$. In the case when $\alpha_r = 1$, big-$M$ has no impact on the equation.

\subsection{Visualization of first- and second-stage solutions of the BAM} \label{visualization_supplement}
\begin{figure*}[h]
  \centering
  \subfloat[\$1000/MWh\label{fig:granular_voll_1000}]{\centering\includegraphics[trim={0.0cm 0cm 4.9cm 0cm},clip, height=1.7in]{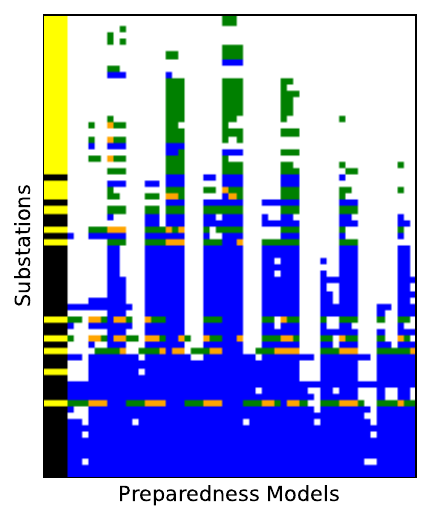}}
  \subfloat[\$2000/MWh\label{fig:granular_voll_2000}]{\centering\includegraphics[trim={0.55cm 0cm 4.9cm 0cm},clip, height=1.7in]{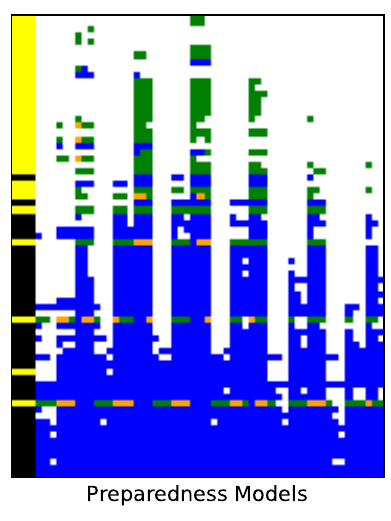}}
  \subfloat[\$3000/MWh\label{fig:granular_voll_3000}]{\centering\includegraphics[trim={0.55cm 0cm 4.9cm 0cm},clip, height=1.7in]{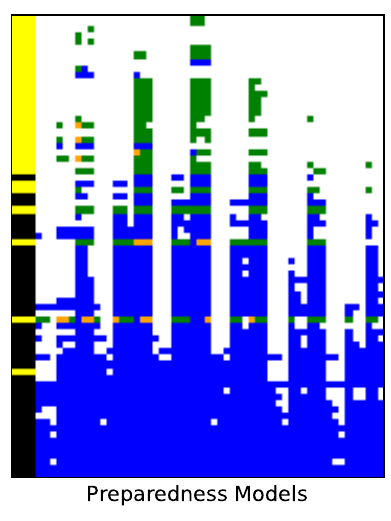}}
  \subfloat[\$4000/MWh\label{fig:granular_voll_4000}]{\centering\includegraphics[trim={0.55cm 0cm 4.9cm 0cm},clip, height=1.7in]{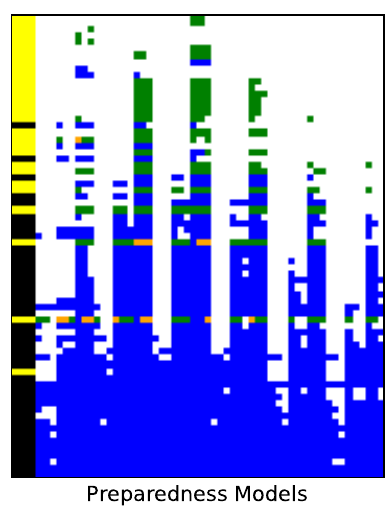}}
  \subfloat[\$5000/MWh\label{fig:granular_voll_5000}]{\centering\includegraphics[trim={0.55cm 0cm 4.9cm 0cm},clip, height=1.7in]{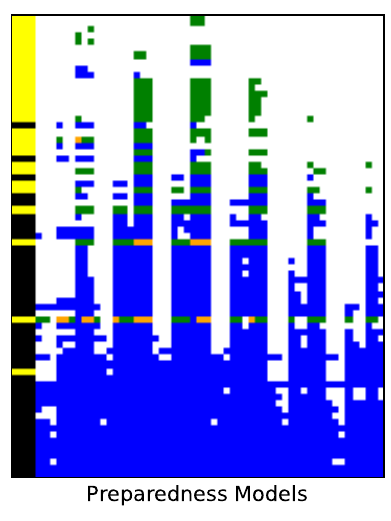}}\hfill
  \subfloat[\$1000/MWh\label{fig:decoupled_1000}]{\centering\includegraphics[trim={0.0cm 0cm 4.9cm 0cm},clip, height=1.7in]{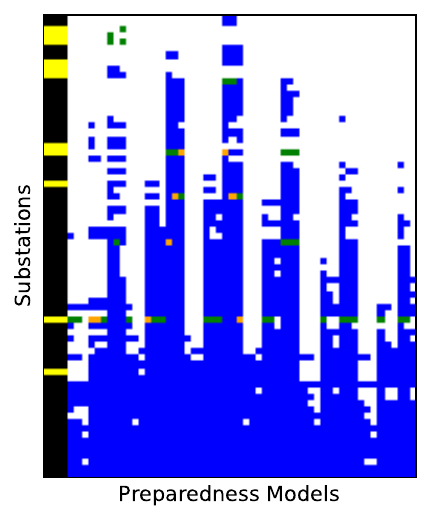}}
  \subfloat[\$2000/MWh\label{fig:decoupled_2000}]{\centering\includegraphics[trim={0.55cm 0cm 4.9cm 0cm},clip, height=1.7in]{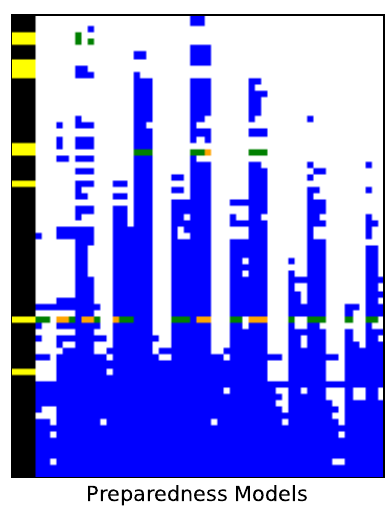}}
  \subfloat[\$3000/MWh\label{fig:decoupled_3000}]{\centering\includegraphics[trim={0.55cm 0cm 4.9cm 0cm},clip, height=1.7in]{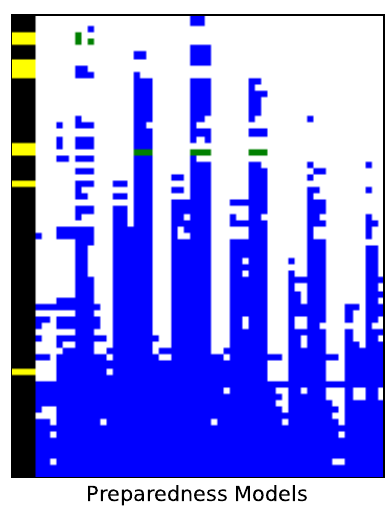}}
  \subfloat[\$4000/MWh\label{fig:decoupled_4000}]{\centering\includegraphics[trim={0.55cm 0cm 4.9cm 0cm},clip, height=1.7in]{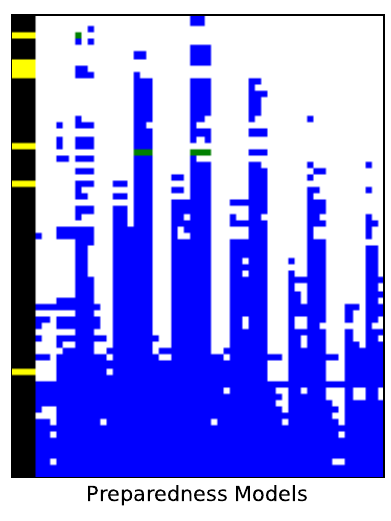}}
  \subfloat[\$5000/MWh\label{fig:decoupled_5000}]{\centering\includegraphics[trim={0.55cm 0cm 4.9cm 0cm},clip, height=1.7in]{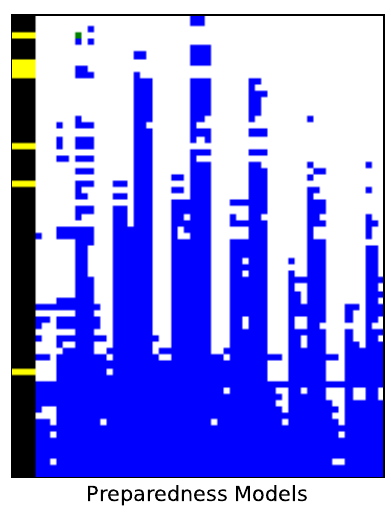}}\hfill
  
  \subfloat[9 feet \label{fig:max_prep_granular_9}]{\centering\includegraphics[trim={0.0cm 0cm 4.9cm 0cm},clip, height=1.7in]{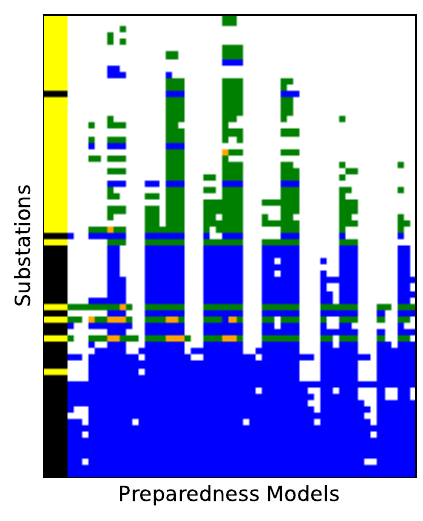}\hspace{0.0cm}}
  \subfloat[12 feet\label{fig:max_prep_granular_12}]{\centering\includegraphics[trim={0.55cm 0cm 4.9cm 0cm},clip, height=1.7in]{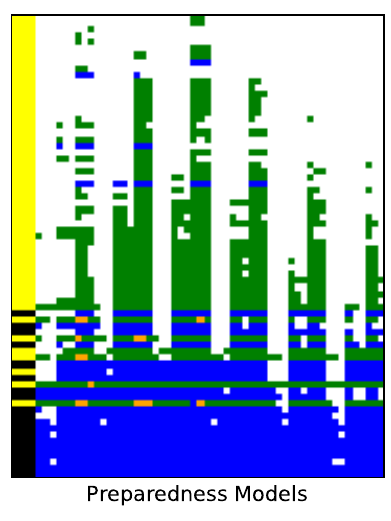}}
  \subfloat[18 feet \label{fig:max_prep_granular_18}]{\centering\includegraphics[trim={0.55cm 0cm 4.9cm 0cm},clip, height=1.7in]{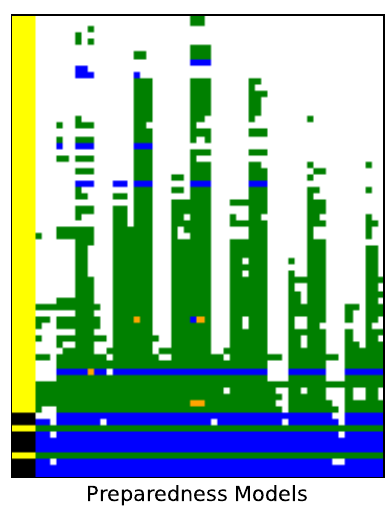}}
  \subfloat[21 feet \label{fig:max_prep_granular_21}]{\centering\includegraphics[trim={0.55cm 0cm 0.0cm 0cm},clip, height=1.7in]{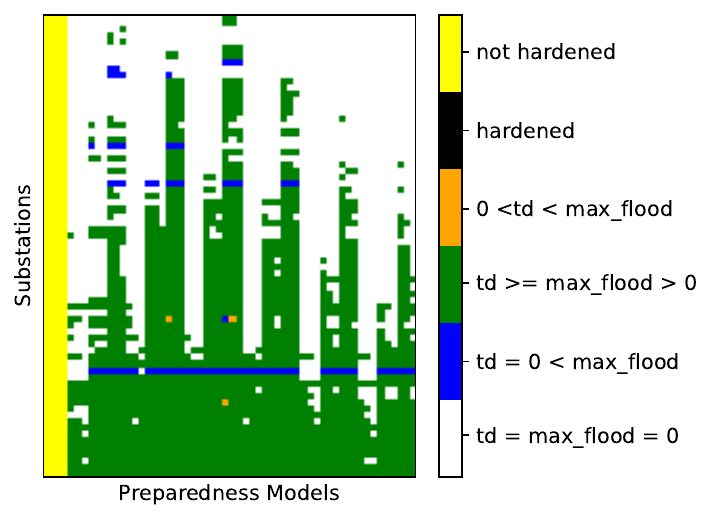}}
  \caption{Plots \eqref{fig:granular_voll_1000} - \eqref{fig:granular_voll_5000} show the optimal first- and second-stage solutions of the BAM for different VOLL values. Plots \eqref{fig:decoupled_1000}-\eqref{fig:decoupled_5000} depict the solutions for different VOLL for decoupled decision-making (modeled as the BAM-D) described in Section \ref{vofc}. Plots \eqref{fig:max_prep_granular_9}-\eqref{fig:max_prep_granular_21} show the solutions for different maximum Tiger Dam effectiveness levels ($\mathcal{T}_i^{\text{max}}$) as discussed in Section \ref{bam_tiger_dam_improvement}. The sub-caption shows the value of the only parameter (VOLL or $\mathcal{T}_i^{\text{max}}$) that is varied from the parameter's default value (\$6,000/MW-hr or 6 ft) to produce the plot. We refer to Section \ref{figure_discussion} for an explanation on how to interpret these figures.}
  \label{fig:overall_label_supp_1}
\end{figure*}

Figure \ref{fig:overall_label_supp_1} presents the first- and second-stage decisions for different parameterizations of the BAM. 

\end{document}